\documentclass{amsart}
\usepackage{amssymb}
\usepackage{microtype}
\usepackage[colorlinks,citecolor=blue,urlcolor=black,linkcolor=black]{hyperref}
\newtheorem{theorem}{Theorem}[section]

\newtheorem{lemma}[theorem]{Lemma}
\newtheorem{corollary}[theorem]{Corollary}
\newtheorem{claim}{Claim}[theorem]
\newtheorem{proposition}[theorem]{Proposition}

\theoremstyle{remark}
\newtheorem{remark}[theorem]{Remark}

\DeclareMathOperator{\PP}{PP}

\DeclareMathOperator{\zfc}{\sf ZFC}

\DeclareMathOperator{\card}{{\mathsf{Card}}}
\DeclareMathOperator{\cf}{cf}

\DeclareMathOperator{\cof}{cof}

\DeclareMathOperator{\cl}{cl}

\DeclareMathOperator{\pcf}{pcf}
\DeclareMathOperator{\tcf}{tcf}

\DeclareMathOperator{\pp}{pp}

\DeclareMathOperator{\Reg}{Reg}
\renewcommand\mid{\mathrel{|}\allowbreak}
\newcommand{\sk}{\vskip.05in}
\title{Normality in Generalized Rudin Spaces}

\date{\today}
\author[T. Eisworth]{Todd Eisworth}
\address{Department of Mathematics, Ohio University, Athens, OH 45701 USA}
\author[M. Kojman]{Menachem Kojman}

\address{Department of Mathematics, Ben-Gurion University of the Negev, P.O.B. 653, Be’er Sheva, 84105 Israel}
\urladdr{https://www.math.bgu.ac.il/~kojman/}
\thanks{The work of the authors was supported by the National Science Foundation under Award No. DMS-2400200 (Todd Eisworth) and by the United States–Israel Binational Science Foundation (BSF) under Grant No. 2023707 (Menachem Kojman).}

\begin{document}

\begin{abstract}
In \cite{Kojman2026}, Kojman generalized Mary Ellen Rudin's construction of a $\zfc$ Dowker space to larger cardinals and gave a pcf-theoretic condition ensuring that the resulting spaces are collectionwise normal. A subsequent manuscript by Eisworth, Kojman, Lambie-Hanson, and Marun \cite{EKLM2026} develops the theory further and leaves open whether these generalized Rudin spaces can fail to be normal. We answer this question by proving that, for cardinals $\mu<\delta$ with $\cf(\delta)=\aleph_0$, the generalized Rudin space $X(\mu,\delta)$ is normal if and only if $\pp(\nu)<\delta$ for every cardinal $\nu\in(\mu,\delta)$ with $\cf(\nu)=\aleph_0$.  Given this characterization, known consistency results provide consistent examples of nonnormal generalized Rudin spaces. 
\end{abstract}
\maketitle

\section{Introduction}

In her 1971 construction~\cite{Rudin1971} of a $\zfc$ Dowker space\footnote{A Dowker space is a normal topological space whose product with the closed unit interval $[0,1]$ is not normal.}, Mary Ellen Rudin considered two related subspaces of the box product $\prod_{0<n<\omega}(\omega_n+1)$:

\[X'=\{x\in {\prod_{0<n<\omega}} (\omega_n+1)\mid (\forall i>0)(\aleph_0<\cf(x(i)))\}\]
 and its subspace
\[
X=\{x\in X'\mid(\exists m)(\forall i>0)(\cf(x(i))<\aleph_m)\}.
\]
She proved that the space $X$ is normal (even collectionwise normal) but not countably paracompact. By Dowker's classical result~\cite{Dowker1951} it follows that $X$ is a normal space for which $X\times[0, 1]$ is nonnormal, that is, a Dowker space.

These two spaces are closely connected. For example, it is easy to see that $X$ is a dense subspace of $X'$ and Rudin herself mentions that $X'$ is the Hewitt realcompactification of $X$.  More important for us, $X'$ possesses strong separation properties that can be exploited to prove things about $X$. This was explicitly noted by Hart~\cite{Hart-notes}, who points out that if we take a discrete family $\mathcal{H}$ of subsets of $X$, then the family $\bar{\mathcal{H}}$ of their closures in $X'$ will be discrete in $X'$. The space $X'$ is paracompact (even ultraparacompact) and hence collectionwise normal, and an expansion of $\bar{\mathcal{H}}$ to disjoint open sets in $X'$ gives us such an expansion of $\mathcal{H}$ in the space $X$. Thus, $X$ inherits its collectionwise normality from $X'$.

Rudin's space remained the only known example in $\zfc$ of a Dowker space for many years, but more examples were discovered in the 1990s: in 1996 Balogh~\cite{Balogh1996} constructed a Dowker space of size $2^{\aleph_0}$, while Kojman and Shelah~\cite{KojmanShelah1998} used the existence of a scale at $\aleph_\omega$ to show that Rudin's space contains a closed Dowker subspace of (absolute) cardinality $\aleph_{\omega+1}$.\footnote{Rudin's space has cardinality $(\aleph_\omega)^{\aleph_0}$, so the sizes of both her example and Balogh's depend on the underlying model of set theory.}

Recently, Kojman~\cite{Kojman2026} generalized Rudin's construction to larger cardinals. This translation replaces $\aleph_0$ by an arbitrary infinite cardinal $\mu$ and replaces the sequence $\langle \aleph_n:n<\omega\rangle$ by a strictly increasing sequence of regular cardinals $\overline{\kappa}=\langle \kappa_n:n<\omega\rangle$ greater than $\mu$. The resulting space $X(\mu,\overline{\kappa})$ turns out to be a Dowker space provided the cardinal $\delta=\sup\overline{\kappa}$ satisfies a certain pcf-theoretic assumption.  This assumption holds automatically when $\delta$ is a singular strong limit cardinal of cofinality  $\omega$, and so this construction answered a question of Rinot by demonstrating that the cardinalities of Dowker spaces form a proper class in $\zfc$.

The study of these generalized Rudin spaces is continued in~\cite{EKLM2026}, where it is shown that the resulting spaces depend only on the choice of $\mu$ and $\delta=\sup\overline{\kappa}$, as different cofinal sequences of regular cardinals above $\mu$ give homeomorphic spaces. Thus, if $\mu<\delta$ are infinite cardinals and $\delta$ has countable cofinality, we can unambiguously define the generalized Rudin space $X(\mu,\delta)$ to be the homeomorphism type of
\begin{equation}\label{defX}
  X(\mu,\overline{\kappa})
  =\left\{x\in\prod_{n<\omega}(\kappa_n+1)
  \mid(\exists m<\omega)(\forall i<\omega)
  (\mu<\cf(x(i))<\kappa_m)\right\}
\end{equation}
with the topology inherited from the box product, where $\overline{\kappa}$ is any $\omega$-sequence (not necessarily increasing or one-to-one) of regular cardinals from the interval $(\mu,\delta)$ with supremum $\delta$. With this notation, Rudin's original Dowker space $X$ is the special case $X(\aleph_0,\aleph_\omega)$, and the main results of~\cite{Kojman2026}, together with the pcf reformulation established below, tell us
\begin{itemize}
  \item $X(\mu,\delta)$ is never countably paracompact, and
    \sk
  \item $X(\mu,\delta)$ is normal (hence a Dowker space) if
    \begin{equation}\label{cond}
      \mu<\nu<\delta\text{ and }\cf(\nu)=\aleph_0
      \Longrightarrow \pp(\nu)<\delta.
    \end{equation}
\end{itemize}

Here ``$\pp$'' refers to Shelah's pseudopower function. For a singular cardinal $\nu$ of countable cofinality, $\pp(\nu)$ is the supremum of the cofinalities $\cf(\prod A/U)$, where $A\subseteq\nu\cap\Reg$ is countable and cofinal in $\nu$ and $U$ is an ultrafilter on $A$ disjoint from the ideal of bounded subsets of $A$. If $\delta$ is a strong limit, then $\pp(\nu)\leq\nu^{\aleph_0}\leq 2^\nu<\delta$ for every $\nu\in(\mu,\delta)$ of countable cofinality, so condition~\eqref{cond} holds.

Condition~\eqref{cond} holds vacuously for every space $X(\mu,\mu^{+\omega})$, since each cardinal strictly between $\mu$ and $\mu^{+\omega}$ is a successor. In particular, Rudin's space $X(\aleph_0,\aleph_\omega)$ satisfies~\eqref{cond} even when $2^{\aleph_0}>\aleph_\omega$. Moreover, assumptions such as {\sf GCH} or the weaker Shelah Strong Hypothesis\footnote{Shelah's Strong Hypothesis {\sf SSH} is that $\pp(\nu)=\nu^+$ whenever $\nu$ is a singular cardinal. This follows from {\sf GCH} and its failure requires consistency of large cardinals.} guarantee that every generalized Rudin space $X(\mu,\delta)$ will be normal and hence Dowker.

The smallest possible counterexample to normality is $X(\aleph_0,\aleph_{\omega+\omega})$, for which condition~(\ref{cond}) reduces to
\[
  \pp(\aleph_\omega)<\aleph_{\omega+\omega}.
\]
By contrast, $X(\aleph_\omega,\aleph_{\omega+\omega})$ satisfies condition~(\ref{cond}) vacuously.

It is therefore natural to ask whether normality must hold in $\zfc$ even when condition~(\ref{cond}) fails, and we are led to the following questions left open by~\cite{Kojman2026} and~\cite{EKLM2026}:

\medskip

\noindent{\bf Main Question}
 {\em Is every Rudin-type space $X(\mu,\delta)$ normal? If not, can the space $X(\aleph_0,\aleph_{\omega+\omega})$ fail to be normal?}

\medskip

In this paper, we answer both of these questions by showing that
condition~(\ref{cond}) is also necessary for $X(\mu,\delta)$ to
be normal. More precisely, we prove the following characterization of normality in terms of pcf theory:

\begin{theorem}\label{thm:normality-main}
  Let $\mu<\delta$ be infinite cardinals with
  $\cf(\delta)=\aleph_0$. Then the following conditions are equivalent:
  \begin{enumerate}
    \item $X(\mu,\delta)$ is normal;
    \sk 
    \item $\max\pcf A<\delta$ for every nonempty countable set
      $A\subseteq\Reg\cap(\mu,\delta)$ with $\sup A<\delta$;
      \sk 
    \item $\pp(\nu)<\delta$ for every cardinal
      $\nu\in(\mu,\delta)$ with $\cf(\nu)=\aleph_0$.
  \end{enumerate}
\end{theorem}
(Note that by~\cite{EKLM2026}, normality is equivalent to collectionwise normality and strong collectionwise normality for these spaces.)

Given the above theorem, we can answer both parts of our main question by quoting known consistency results to show that the smallest previously unresolved case can indeed fail to be normal. 

\begin{corollary}\label{cor:nonnormal-consistency}
  It is consistent, relative to suitable large cardinal assumptions,
  that the space $X(\aleph_0,\aleph_{\omega+\omega})$ is not normal.
\end{corollary}

Any model of $\zfc$ in which $\aleph_\omega$ is strong limit and
$2^{\aleph_\omega}>\aleph_{\omega+\omega}$ satisfies
$\pp(\aleph_\omega)>\aleph_{\omega+\omega}$, so
$X(\aleph_0,\aleph_{\omega+\omega})$ is not normal by
Theorem~\ref{thm:normality-main}. In 1983 Shelah~\cite[Theorem~2.2(1)]{Shelah1983}
constructed such models from a supercompact cardinal. Gitik and
Magidor~\cite{GitikMagidor1992,GitikMagidor1994} later gave a different
construction using extender-based forcing and also arranging for GCH to hold below $\aleph_\omega$, something that Shelah's construction did not achieve.

The strong limit assumption is not necessary for the failure of normality, however, because adding Cohen reals to such a model will preserve
the inequality
$\pp(\aleph_\omega)>\aleph_{\omega+\omega}$.
This shows that the space $X(\aleph_0,\aleph_{\omega+\omega})$ can fail to be normal
in models in which the continuum is arbitrarily large.

Finally, the pcf criterion yields topological characterizations of
normality in terms of cardinal invariants. In
Section~\ref{sec:rudin-cardinal-functions}, we compute the global
pseudocharacter $\psi(X)$, character $\chi(X)$, and tightness $t(X)$
of these spaces. Together with Theorem~\ref{thm:normality-main},
these calculations give the following corollary.

\begin{corollary}\label{cor:topological-intro}
  Let $\mu<\delta$ be infinite cardinals with
  $\cf(\delta)=\aleph_0$, and put $X=X(\mu,\delta)$.
  Then the following conditions are equivalent:
  \begin{enumerate}
    \item $X$ is normal;
    \sk
    \item $\psi(X)=\chi(X)$;
    \sk
    \item $\psi(X)=t(X)$;
    \sk
    \item $\psi(X)=\chi(X)=t(X)=\delta$.
  \end{enumerate}
\end{corollary}

In fact, $\psi(X)=\delta$ and $\chi(X)=t(X)$ always hold for these
spaces. Thus failure of normality is witnessed by the strict inequality
$\psi(X)<\chi(X)=t(X)$.

\section{Notation and preliminaries}\label{sec:background}

Our notation is standard, so in this section we will only sketch the concepts we need and give references to the literature.

\subsection{Products and Rudin-type spaces}

We order ordinal-valued functions pointwise, with $f<g$ meaning
$f(i)<g(i)$ at every coordinate, and define

$$(h,x]=\{y:h<y\leq x\}.$$
For a set $A$ of regular cardinals, write
$\prod A=\prod_{\kappa\in A}\kappa$ and let $\cf(\prod A,<)$
denote the least size of a pointwise cofinal subset.

Let $\mu$ be an infinite cardinal and
$\overline\kappa=\langle\kappa_n:n<\omega\rangle$ a sequence of
regular cardinals above $\mu$ with no largest term. The terms may repeat
and may occur in any order. Put $\delta=\sup_n\kappa_n$.
Define $X(\mu,\overline\kappa)$ by~\eqref{defX} and put
\[
 X'(\mu,\overline\kappa)=
 \left\{x\in\prod_{n<\omega}(\kappa_n+1):
 (\forall n<\omega)\,\mu<\cf(x(n))\right\}.
\]
Both spaces inherit the box topology from the order topologies.
In either space, intersecting $(h,x]$ with the space, for $h<x$,
gives a neighborhood base at $x$.

For a countably infinite set $C$ of regular cardinals above $\mu$
with no maximum, $X(\mu,C)$ and $X'(\mu,C)$ denote the corresponding
$C$-indexed spaces.
By~\cite{EKLM2026}, the homeomorphism type of the pair
$(X'(\mu,\overline\kappa),X(\mu,\overline\kappa))$ depends only on
$\mu$ and $\delta$, giving the notation $X(\mu,\delta)$.
\subsection{Topological conventions}

For $y\in Y$, let $\chi(y,Y)$ denote its character and $\psi(y,Y)$ its
pseudocharacter. Thus $\chi(y,Y)$ is the least size of a neighborhood base
at $y$, and $\psi(y,Y)$ is the least size of a family of open neighborhoods
whose intersection is $\{y\}$. The tightness $t(y,Y)$ is the least infinite
cardinal $\kappa$ such that, whenever $y\in\overline D$ for $D\subseteq Y$,
there is $E\subseteq D$ with $|E|\leq\kappa$ and $y\in\overline E$.
We take these invariants to be infinite and put
\[
 \chi(Y)=\sup_{y\in Y}\chi(y,Y),
 \qquad
 \psi(Y)=\sup_{y\in Y}\psi(y,Y),
 \qquad
 t(Y)=\sup_{y\in Y}t(y,Y).
\]

\subsection{Pcf preliminaries}

Let $J$ be a proper ideal on a set $A$ of regular cardinals. For
$f,g\in\prod A$, define
\[
\begin{split}
 f<_Jg &\quad\Longleftrightarrow\quad
 \{\kappa\in A:f(\kappa)\geq g(\kappa)\}\in J,\\
 f=_Jg &\quad\Longleftrightarrow\quad
 \{\kappa\in A:f(\kappa)\neq g(\kappa)\}\in J.
\end{split}
\]
For the ideal $\mathrm{Fin}$ of finite subsets, we write $<^*$ and $=^*$.
The ideal of subsets of $A$ bounded in $\sup A$ is denoted by
$J^{\mathrm{bd}}[A]$. If $A$ has order type $\omega$, then
$J^{\mathrm{bd}}[A]=\mathrm{Fin}$.

If $\lambda$ is regular and $\prod A$ has a $<_J$-increasing,
$<_J$-cofinal sequence $\langle f_\alpha:\alpha<\lambda\rangle$, we write
\[
 \tcf(\prod A/J)=\lambda
\]
and call the sequence a $\lambda$-\emph{scale}.
The true cofinality need not exist. For an ultrafilter $U$, the notation
$\prod A/U$ uses the ideal dual to $U$.

A nonempty set $A\subseteq\Reg$ is \emph{progressive} if $|A|<\min A$.
For such a set, put
\[
 \pcf(A)=\{\cf(\prod A/U):U\text{ is an ultrafilter on }A\}.
\]
The pcf theorem~\cite[Proposition~2.19]{Gilton2022} gives a largest
member of $\pcf(A)$ and the identity
\begin{equation}\label{eq:pcf-pointwise}
 \cf(\prod A,<)=\max\pcf(A).
\end{equation}
Also, $A\subseteq B$ implies $\pcf(A)\subseteq\pcf(B)$.
Our applications concern nonempty countable sets of uncountable regular cardinals,
which are progressive. We refer to~\cite{Shelah1994,Kojman2001} for these
facts.

We use pseudopowers only at singular cardinals of countable cofinality.
For such a cardinal $\nu$, let $\operatorname{PP}(\nu)$ consist of the
regular cardinals of the form $\tcf(\prod A/J)$, where
$A\subseteq\nu\cap\Reg$ is countable and cofinal in $\nu$,
and $J$ is a proper ideal on $A$ containing $J^{\mathrm{bd}}[A]$. Then
\[
 \pp(\nu)=\sup\operatorname{PP}(\nu).
\]
We use the following facts; see~\cite[Chapter~II]{Shelah1994}.
\begin{enumerate}
\item There is a cofinal set $A\subseteq\nu\cap\Reg$ of order type
$\omega$ such that $\tcf(\prod A/\mathrm{Fin})=\nu^+$.
In particular, $\nu^+\in\operatorname{PP}(\nu)$.

\item If $\lambda\in\operatorname{PP}(\nu)$ and
$\nu<\theta<\lambda$ is regular, then
$\theta\in\operatorname{PP}(\nu)$ (No Holes).

\item If $\theta<\eta\leq\pp(\theta)$ and
$\cf(\theta)=\cf(\eta)=\aleph_0$, then inverse
monotonicity~\cite[Proposition~1.4.3]{Eisworth2023} gives
$\PP(\eta)\subseteq\PP(\theta)$ and hence $\pp(\eta)\leq\pp(\theta)$.
\end{enumerate}

We finish with the normality criterion from~\cite{Kojman2026}.
\begin{theorem}[Kojman]\label{thm:kojman-normality}
Let $\mu$ be an infinite cardinal and let
$\overline\kappa=\langle\kappa_n:n<\omega\rangle$ be a strictly increasing
sequence of regular cardinals above $\mu$. Suppose that, for every
$n<\omega$ and every nonempty countable
$A\subseteq\Reg\cap(\mu,\kappa_n)$,
\[
 \cf\left(\prod A,<\right)\leq\kappa_n.
\]
Then $X(\mu,\overline\kappa)$ is collectionwise normal.
\end{theorem}
By~\eqref{eq:pcf-pointwise}, the hypothesis can also be written
$\max\pcf(A)\leq\kappa_n$.

\section{The pcf characterization of normality}

This section will establish Theorem~\ref{thm:normality-main}.  We start with some easy lemmas that show us that the pcf-theoretic content of the theorem can be formulated in a few equivalent ways.   

\begin{lemma}\label{lem:pcf-local}
Let $\mu<\gamma$ be infinite cardinals. The following are equivalent.
\begin{enumerate}
\sk
\item $\max\pcf(A)\leq\gamma$ for every countable set
$A\subseteq(\mu,\gamma)\cap\Reg$ with $\sup A<\gamma$.
\sk
\item $\pp(\nu)\leq\gamma$ for every cardinal $\nu\in(\mu,\gamma)$
of countable cofinality.
\end{enumerate}
\end{lemma}

\begin{proof}
Suppose $A$ witnesses the failure of~\textup{(1)}. Put
$\lambda=\max\pcf(A)$ and choose an ultrafilter $U$ on $A$ witnessing
$\cf(\prod A/U)=\lambda$. Since $\lambda>\gamma>\sup A$, $U$ is
nonprincipal. Let $\nu$ be the least ordinal such that $A\cap\nu\in U$.
Then $\nu$ is a limit point of $A$, so it is a cardinal of countable
cofinality with $\mu<\nu\leq\sup A<\gamma$. By minimality, the
restriction of $U$ to $A\cap\nu$ contains no bounded subset of $\nu$.
It still witnesses cofinality $\lambda$, so $\lambda\in\PP(\nu)$ and
$\pp(\nu)\geq\lambda>\gamma$.

Conversely, suppose $\nu$ witnesses the failure of~\textup{(2)}.
Choose $\lambda\in\PP(\nu)$ above $\gamma$ and let $(A,J)$ witness
this representation. Since $J$ contains the bounded ideal, we may
discard the coordinates at or below $\mu$ and assume that
$A\subseteq(\mu,\nu)\cap\Reg$. Extending the dual filter of $J$ to an
ultrafilter gives $\lambda\in\pcf(A)$. Thus $\sup A=\nu<\gamma$ and
$\max\pcf(A)\geq\lambda>\gamma$.
\end{proof}

Moving on to the statement of Theorem~\ref{thm:normality-main}, we make use now of the countable cofinality of~$\delta$.   The first two clauses
of the following lemma express the negations of items~(2) and~(3) of
Theorem~\ref{thm:normality-main}, respectively, since $\max\pcf(A_0)$
is regular and $\delta$ is singular. We will use the two scales at the end of the lemma to prove nonnormality.

\begin{lemma}\label{lem:pcf}\label{lem:pcf-nonnormality}
Suppose $\mu<\delta$ are infinite cardinals and
$\operatorname{cf}(\delta)=\aleph_0$. The following are equivalent.
\begin{enumerate}
\sk
\item There is a countable set
$A_0\subseteq(\mu,\delta)\cap\operatorname{Reg}$ such that
\[
    \sup A_0<\delta<\max\operatorname{pcf}(A_0).
\]
\sk
\item There is a cardinal $\nu\in(\mu,\delta)$ of countable
cofinality such that $\operatorname{pp}(\nu)\geq\delta$.

\sk
\item There is a cardinal $\nu\in(\mu,\delta)$ of countable
cofinality such that $\delta^+\in\operatorname{PP}(\nu)$.

\sk
\item\label{item:pcf-scales} There are countable sets $A$ and $B$ of regular cardinals
and a proper ideal $J$ on $A$ extending $J^{\mathrm{bd}}[A]$
such that
\[
    \mu^+<\min A,\qquad
    \sup A<\min B,\qquad
    \sup B=\delta,\qquad
    \operatorname{otp}(B)=\omega,
\]
and
\[
    \operatorname{tcf}\bigl(\prod A/J\bigr)
    =\operatorname{tcf}\bigl(\prod B/\mathrm{Fin}\bigr)
    =\delta^+.
\]
The corresponding scales can be chosen so that each function assumes only values of cofinality $\mu^+$.
\end{enumerate}
\end{lemma}

\begin{proof}
By Shelah's scale theorem~\cite[Chapter~II]{Shelah1994}, fix a
cofinal set $B_0\subseteq\delta\cap\operatorname{Reg}$ of order
type $\omega$ such that
\[
    \operatorname{tcf}\bigl(\prod B_0/\mathrm{Fin}\bigr)=\delta^+.
\]

\noindent $(1)\Rightarrow(2)$. Choose an ultrafilter $U$ on $A_0$
witnessing $\lambda=\max\pcf(A_0)$, and let $\nu$ be the least ordinal
with $A_0\cap\nu\in U$. As in the proof of Lemma~\ref{lem:pcf-local},
$\mu<\nu\leq\sup A_0$, $\cf(\nu)=\aleph_0$, and the restriction of $U$
to $A_0\cap\nu$ witnesses $\lambda\in\PP(\nu)$. Thus
$\pp(\nu)\geq\lambda>\delta$, and $A_0\cap\nu$ witnesses~\textup{(1)}
with supremum $\nu$.

\medskip

\noindent $(2)\Rightarrow(3)$. Since $\nu<\delta\leq\operatorname{pp}(\nu)$
and both $\nu$ and $\delta$ have countable cofinality, inverse
monotonicity~\cite[Proposition~1.4.3]{Eisworth2023} gives
\[
    \operatorname{PP}(\delta)\subseteq\operatorname{PP}(\nu).
\]
The scale on $B_0$ witnesses $\delta^+\in\operatorname{PP}(\delta)$,
so the same $\nu$ witnesses~\textup{(3)}.

\medskip

\noindent $(3)\Rightarrow(4)$. Choose a countable cofinal set
$A\subseteq\nu\cap\operatorname{Reg}$ and a proper ideal
$J\supseteq J^{\mathrm{bd}}[A]$ witnessing
$\operatorname{tcf}(\prod A/J)=\delta^+$. Discarding the bounded
initial segment of $A$ at or below $\mu^+$ and restricting $J$
accordingly, we may assume $\mu^+<\min A$. Let
$B=\{b\in B_0:b>\nu\}$. This removes only finitely many terms
from $B_0$, so the restricted scale still has length $\delta^+$,
and $\sup A=\nu<\min B$.

For each $\kappa\in A\cup B$ and $\xi<\kappa$, put
\[
    e_\kappa(\xi)=\mu^+\cdot(\xi+1),
\]
where the product is ordinal multiplication. Since $\kappa>\mu^+$,
$e_\kappa$ is a strictly increasing map from $\kappa$ to $\kappa$
with $e_\kappa(\xi)>\xi$ and
$\operatorname{cf}(e_\kappa(\xi))=\mu^+$. Replacing each scale
function $f$ by the function $\kappa\mapsto e_\kappa(f(\kappa))$
therefore preserves strict increase modulo the respective ideal
and preserves cofinality.

\medskip

\noindent $(4)\Rightarrow(1)$. Any ultrafilter on $A$ disjoint from $J$
preserves the scale on $A$. Thus $\delta^+\in\operatorname{pcf}(A)$,
and $A_0=A$ witnesses~\textup{(1)}.
\end{proof}

The preceding lemma tells us that the two pcf conditions in Theorem~\ref{thm:normality-main} are equivalent, so what remains is to prove that these give necessary and sufficient conditions for the associated generalized Rudin space to be normal. 

\subsection*{{\sf Sufficiency:}}

Assume (3) of Theorem~\ref{thm:normality-main}, so $\pp(\nu)<\delta$
for every cardinal $\nu\in(\mu,\delta)$ of countable cofinality. By
Lemma~\ref{lem:pcf-local}, a regular $\kappa\in(\mu,\delta)$ satisfies
Kojman's bound precisely when $\pp(\eta)\leq\kappa$ for every cardinal
$\eta\in(\mu,\kappa)$ of countable cofinality. We show that such
$\kappa$ are cofinal in $\delta$.

Fix a regular cardinal $\gamma\in(\mu,\delta)$. If $\pp(\eta)\leq\gamma$
for every cardinal $\eta\in(\mu,\gamma)$ of countable cofinality, take
$\kappa=\gamma$. Otherwise, let $\theta$ be the least cardinal of
countable cofinality in $(\mu,\gamma)$ with $\pp(\theta)>\gamma$, and put
\[
    \tau=\pp(\theta),\qquad \kappa=\tau^+.
\]
By (3), $\tau<\delta$, and hence $\kappa<\delta$ since $\delta$ is a
limit cardinal. For $\eta\in(\mu,\theta)$ of countable cofinality,
minimality gives $\pp(\eta)\leq\gamma$. For $\theta<\eta\leq\tau$ of
countable cofinality, inverse monotonicity gives $\pp(\eta)\leq\tau$.
Thus $\kappa$ works in this case as well.

As $\gamma$ was arbitrary, we can choose a strictly increasing sequence
$\overline\kappa=\langle\kappa_n:n<\omega\rangle$ of such regular
cardinals cofinal in $\delta$. By Theorem~\ref{thm:kojman-normality},
$X(\mu,\overline\kappa)$ is normal, and its homeomorphism
type is $X(\mu,\delta)$, as required.

\subsection*{{\sf Necessity:}}

We prove $(1)\Rightarrow(3)$ of Theorem~\ref{thm:normality-main} by
contraposition. The model for our proof is the classic example (due independently to  Dieudonn\'e and Morse)
of the product of ordinal spaces $(\omega_1+1)\times\omega_1$. In their setting, the vertical edge
$\{\omega_1\}\times\omega_1$ and the diagonal are disjoint closed
sets. An application of Fodor's lemma then shows that every open neighborhood of the diagonal
contains a final square $(\beta,\omega_1)^2$, and  that every open
neighborhood of the edge meets this square. In our setting,  we will form a
diagonal by pairing points with the same index on the two scales from
Lemma~\ref{lem:pcf}, and then interleave two sequences with a common
pointwise supremum, as in Kelley's treatment of the ordinal
example~\cite[Chapter~4, Problem~E]{Kelley1975}. For related pressing-down
arguments in products, see~\cite{KojmanShelah1998,KojmanMichalewski}.

\begin{proposition}
If $\pp(\nu)\geq\delta$ for some cardinal $\nu\in(\mu,\delta)$ of
countable cofinality, then $X(\mu,\delta)$ is not normal.
\end{proposition}

\begin{proof}
Put $\lambda=\delta^+$. By Lemma~\ref{lem:pcf}\textup{(\ref{item:pcf-scales})}, choose $A,B,J$ and scales
\[
    \langle p_\alpha:\alpha<\lambda\rangle\subseteq\prod A,
    \qquad \langle q_\alpha:\alpha<\lambda\rangle\subseteq\prod B
\]
whose coordinate values all have cofinality $\mu^+$. The first scale
is increasing and cofinal in $\prod A$ modulo $J$, and the second is strictly increasing and cofinal in $\prod B$ modulo the finite ideal.

Put $C=A\cup B$ and use $X=X(\mu,C)$ as our presentation of
$X(\mu,\delta)$.  Let $\Delta(c)=c$ be the diagonal function for $c\in C$. Then
$\Delta\in X'(\mu,C)\setminus X$, since the cofinalities of its
coordinates are unbounded in $\delta$.

Our analogue of the vertical edge is the closed set
\begin{equation}
    \label{eqn:Fdef}
    F=\{x\in X:(\forall a\in A)\,x(a)=a\}.
\end{equation}
Indeed, $F$ is the intersection, over $a\in A$, of the preimages of the closed singletons
$\{a\}\subseteq a+1$ under the continuous coordinate projections $x\mapsto x(a)$.
Replacing the $A$-coordinates of any $x\in X$ by their top values
leaves their cofinalities bounded by $\sup A<\delta$, so $F$ is
cofinal in $(X,\leq)$. On the other hand, every $x\in X$ satisfies
$x(b)<b$ for all but finitely many $b\in B$, since $B$ has order
type $\omega$ and is cofinal in $\delta$.

For our analogue of the diagonal, we put $r_\alpha=p_\alpha\cup q_\alpha$, and define

\begin{equation}
    \label{eqn:Ddef}
 D=\bigcup_{\alpha<\lambda}\left\{z\in\prod C:
   \begin{array}{l}
     \cf(z(c))=\mu^+\text{ for every }c\in C,\\
     z\restriction A=_Jp_\alpha\text{ and }z\restriction B=^*q_\alpha
   \end{array}\right\}.
\end{equation}
Thus, membership in the set $D$ is unchanged if we make a small modification to an element, provided we keep control of the cofinalities: on the coordinates in $A$ we are allowed to modify $r_\alpha$ on a $J$-small set, while we are allowed to make finite changes on the coordinates from $B$.  
The common cofinality bound gives $D\subseteq X$, and we define

\begin{equation}
    \label{eqn:Gdef}
    G = \cl_X D.
\end{equation}
The closed subsets $F$ and $G$ of $X$ are intended to refute the normality of the space. The next claim tells us that these are indeed disjoint closed subsets of $X$:

\begin{claim}
$F\cap G=\emptyset$.
\end{claim}

\begin{proof}
Fix $x\in F$. Since $x(b)<b$ for all but finitely many $b\in B$,
choose $\alpha<\lambda$ such that $x\restriction B<^*q_\alpha$.
Put $h\restriction A=p_\alpha$ and let $h\restriction B$ be
identically zero. Then $h<x$, and we claim that $(h,x]\cap X$ misses
$D$. Given $z\in D$, choose $\beta$ witnessing its membership. If
$\beta\leq\alpha$, then $z(a)\leq p_\alpha(a)=h(a)$ outside a
member of $J$, hence at some $a\in A$. If $\alpha<\beta$, then
\[
    x\restriction B<^*q_\alpha<^*q_\beta=^*z\restriction B,
\]
so $z(b)>x(b)$ for some $b\in B$. In either case $z\notin(h,x]$.
Thus $x\notin G$.
\end{proof}

Our next step is to show that every open neighborhood of the closed diagonal $G$ contains an end segment of $X\cap\prod C$.

\begin{claim}\label{endseg}
If $U$ is an open neighborhood of $G$, then there is a function $h\in\prod C$
such that $\{z\in X\cap\prod C:h<z\}\subseteq U$, that is, such that  
\begin{equation}
    X\cap\prod_{c\in C}(h(c), c)\subseteq U.
\end{equation}
\end{claim}

\begin{proof}
First note that $D$ is pointwise cofinal in $\prod C$. Indeed, given  some $t\in\prod C$,
we can choose $\alpha<\delta^+$ so that $t\restriction A<_Jp_\alpha$ and
$t\restriction B<^*q_\alpha$. Wherever $r_\alpha(c)\leq t(c)$,
replace $r_\alpha(c)$ by $t(c)+\mu^+$. This ordinal is below $c$ and
has cofinality $\mu^+$. The resulting function belongs to $D$ and is
everywhere above $t$.

Now suppose the claim fails for some open neighborhood $U$ of $G$. Then
$$E=\left(X\cap\prod C\right)\setminus U$$ is also pointwise cofinal in $\prod C$.
Since each $c\in C$ is regular and greater than $\mu^+$, we can
recursively choose $d_\alpha\in D$ and $e_\alpha\in E$ for
$\alpha<\mu^+$ so that
\[
    d_\alpha<e_\alpha<d_\beta\qquad(\alpha<\beta<\mu^+).
\]
Their common pointwise supremum $x$ satisfies $x(c)<c$ and
$\cf(x(c))=\mu^+$ for every $c\in C$, so $x\in X$. For each $h<x$,
the countability of $C$ and regularity of $\mu^+$ give some $\alpha<\mu^+$
with $h<d_\alpha<e_\alpha\leq x$. Thus $x\in\cl_X D=G\subseteq U$,
while $x\in\cl_X E\subseteq X\setminus U$, a contradiction.
\end{proof}

Finally, we show that every open neighborhood of our set $F$ meets every end segment of $X\cap\prod C$.  

\begin{claim}
\label{claim:fat}
Every open neighborhood of $F$ meets $\{z\in X\cap\prod C:h<z\}$
for every $h\in\prod C$.
\end{claim}

\begin{proof}
Let $V$ be an open neighborhood of $F$ and fix $h\in\prod C$.
Define $y\in F$ by $y(a)=a$ for $a\in A$ and $y(b)=h(b)+\mu^+$
for $b\in B$. Choose $g<y$ with $(g,y]\cap X\subseteq V$ and set
\[
    z(a)=\max\{g(a),h(a)\}+\mu^+\quad(a\in A),\qquad
    z(b)=y(b)\quad(b\in B).
\]
Every coordinate of $z$ has cofinality $\mu^+$, so $z\in X\cap\prod C$.
Since $h<z$ and $g<z\leq y$, this $z$ belongs to $V$ and to the
required tail.
\end{proof}

The conjunction of Claims~\ref{endseg} and~\ref{claim:fat} tells us that every open neighborhood of $F$ must meet every open neighborhood of $G$, and therefore the disjoint closed sets $F$ and $G$ witness that $X$ is not normal.
\end{proof}

This proves $(1)\Rightarrow(3)$ and completes the proof of
Theorem~\ref{thm:normality-main}. In the construction above, both $F$
and $G$ have $\Delta$ in their closure in $X'(\mu,C)$, as witnessed
by the points $y$ in the last claim and by the cofinality of $D$,
respectively. The open neighborhood $X\setminus G$ of $F$ misses
$D$, so it contains no end segment of $X\cap\prod C$. For comparison
with the end-segment property in Rudin's original space,
see~\cite[Lemma~2.25]{KojmanMichalewski}.

\section{Topological characterization of normality}
\label{sec:rudin-cardinal-functions}

In the preceding section we showed that \(X(\mu,\delta)\) is normal precisely when \(\max\operatorname{pcf}(A)<\delta\) for every nonempty bounded countable set \(A\) of regular cardinals in \((\mu,\delta)\). We now show how this condition has a direct translation in terms of cardinal functions defined from the topology of $X(\mu,\delta)$.
For each such \(A\), there is a point whose pseudocharacter is \(\sup A\) and whose character and tightness are both \(\max\operatorname{pcf}(A)\). Passing to the whole space, we find that its pseudocharacter is always~\(\delta\), while its character and tightness exceed \(\delta\) precisely when the pcf condition fails.

Throughout this section, $\mu<\delta$ are infinite cardinals,
$\cf(\delta)=\omega$, and we assume that $X=X(\mu,\delta)$ has the presentation
$X(\mu,\overline\kappa)$ where $\overline{\kappa}$ is strictly increasing.
In this context, we define
\[
 \begin{split}
 \mathcal S_{\mu,\delta}
   &=\{\nu:\nu\in (\mu,\delta)\cap\card\text{ and } \cf(\nu)=\omega<\nu\},\\
   \\
 H_{\mu,\delta}&=\sup\{\pp(\nu):\nu\in\mathcal S_{\mu,\delta}\},
 \end{split}
\]
where we set $\sup\varnothing=0$. Our next goal is the following computation of the character, pseudocharacter, and tightness of $X= X(\mu,\delta)$. 

\begin{theorem}\label{thm:global}
The space $X=X(\mu,\delta)$ satisfies
\begin{enumerate}
    \item $\psi(X)=\delta$, and
    \sk 
    \item $\chi(X)=t(X)=\max\{\delta, H_{\mu,\delta}\}$.
\end{enumerate}
\end{theorem}

The term $\delta$ in the second clause is necessary in general.
For example, if $\delta=\mu^{+\omega}$, then
$\mathcal S_{\mu,\delta}=\varnothing$, while all three global
invariants equal $\delta$. The next lemma computes all three invariants at a nonisolated point.

\begin{lemma}\label{lem:local}
For $n<\omega$, define
$T_n=\{\alpha\leq\kappa_n:\cf(\alpha)>\mu\}$.
If $x\in X$ is not an isolated point, we define
\begin{equation*}
    I_x := \{n<\omega: x(n)\text{ is nonisolated in }T_n\},
\end{equation*}
and
\begin{equation*}
    A_x:= \{\cf(x(n)):n\in I_x\}.
\end{equation*}
Then $A_x$ is a nonempty countable set of regular cardinals that is bounded below $\delta$, and
\begin{equation}\label{eq:local}
 \psi(x,X)=\sup A_x,
 \end{equation}
while
 \begin{equation}
 \chi(x,X)=t(x,X)=\max\pcf(A_x).
\end{equation}
\end{lemma}

\begin{proof}
Since $x$ is nonisolated, we know $I_x\ne\varnothing$, and
$\sup A_x<\delta$ by~\eqref{defX}. For each $n\in I_x$ choose a
cofinal $D_n\subseteq x(n)$ of order type $\cf(x(n))$ consisting of
ordinals of cofinality~$\mu^+$. Indeed, for $\alpha<x(n)$ the first
ordinal of cofinality greater than $\mu$ above $\alpha$ is
$\alpha+\mu^+$, and nonisolation puts it below $x(n)$.

Put $D=\prod_{n\in I_x}D_n$, viewing each member $d\in D$ as a point of $X$ by
setting $d(n)=x(n)$ for every $n\notin I_x$. For $d\in D$, let $U_d$ consist of those $y\in X$ that agree with $x$
off $I_x$ and satisfy $d(n)<y(n)\leq x(n)$ for every $n\in I_x$.
Since $x(n)$ is isolated in $T_n$ off $I_x$ and each $D_n$ is cofinal
in $x(n)$ on $I_x$, the boxes $U_d$ form a neighborhood base at $x$. A family of these boxes isolates $x$ exactly when its lower  bounds are cofinal in each $x(n)$ for $n\in I_x$. If they are bounded in
one coordinate, replacing $x(n)$ by a member of $D_n$ above the bound
gives a second point in the intersection. Since $I_x$ is countable,
$\psi(x,X)=\sup A_x$.

For $E\subseteq D$, the same boxes give $x\in\cl_X E$ exactly when $E$
is pointwise cofinal in $D$. Hence
\[
 \cof(D,\leq)\leq t(x,X)\leq\chi(x,X)\leq\cof(D,\leq).
\]
Identify each $D_n$ with its order type. The functions constant on
coordinates of equal cofinality form a cofinal copy of $\prod A_x$ in
$D$, since every countable subset of $\theta\in A_x$ is bounded in
$\theta$. Thus $\cof(D,\leq)=\max\pcf(A_x)$ by~\eqref{eq:pcf-pointwise}.
\end{proof}

\begin{proof}[Proof of Theorem~\ref{thm:global}]
Given a nonempty countable $A\subseteq(\mu,\delta)\cap\Reg$ with
$\sup A<\delta$, we first find a point $x\in X$ with $A_x=A$.
Fix an $\omega$-sequence $\langle\theta_i:i<\omega\rangle$ with
range~$A$, allowing repetitions. In our fixed increasing presentation
$X=X(\mu,\bar\kappa)$, let
\[
    n^*=\min\{n<\omega:\sup A<\kappa_n\}.
\]
Set $x(n)=\mu^+$ for $n<n^*$ and
\[
    x(n^*+i)=\mu^+\cdot\theta_i \qquad(i<\omega),
\]
where the product is ordinal multiplication. The ordinal
$\mu^+\cdot\theta_i$ has cardinality and cofinality $\theta_i$,
so it lies below $\kappa_{n^*+i}$, and $x$ is in $X$ because the cofinalities
of its values are at most $\sup A$, hence bounded below $\delta$. 

For each $i<\omega$, the ordinals $\mu^+\cdot(\xi+1)$ for
$\xi<\theta_i$ have cofinality $\mu^+$ and form a sequence in
$T_{n^*+i}$ increasing to $x(n^*+i)$, so $x(n^*+i)$ is
nonisolated in $T_{n^*+i}$ and hence $n^*+i\in I_x$.  
Since $x(n)=\mu^+$ is isolated in $T_n$ for $n<n^*$, we have
$I_x=\{n^*+i:i<\omega\}$. Hence
\[
    A_x=\{\cf(x(n^*+i)):i<\omega\}
       =\{\theta_i:i<\omega\}=A.
\]

By Lemma~\ref{lem:local}, every nonisolated point of $X$ has
pseudocharacter below $\delta$. Taking $A=\{\theta\}$ gives a point
of pseudocharacter $\theta$ for every regular $\theta\in(\mu,\delta)$.
Since these cardinals are cofinal in $\delta$, we have $\psi(X)=\delta$.
Taking suprema of the pointwise character and tightness gives
\begin{equation}\label{eq:bounded}
 \chi(X)=t(X)=
 \sup\left\{\max\pcf(A):
 \begin{array}{l}
 \varnothing\neq A\subseteq\Reg\cap(\mu,\delta),\\
 |A|\leq\aleph_0,\quad \sup A<\delta
 \end{array}\right\}.
\end{equation}
In particular, this supremum is at least $\delta$.

We can now compute \eqref{eq:bounded} in terms of pseudopowers.
Let $A$ occur in that supremum and let
$\lambda=\max\pcf(A)$. Choose an ultrafilter $U$ on $A$ such that
$\cf(\prod A/U)=\lambda$. If $U$ is principal, then
$\lambda\in A$, so $\lambda<\delta$.
Otherwise, let $\nu$ be the least ordinal for which $A\cap\nu\in U$,
and put $B=A\cap\nu$. If $\nu=\beta+1$, then nonprincipality gives
$A\cap\beta=B\setminus\{\beta\}\in U$, contrary to the choice of $\nu$.
Thus $\nu$ is a limit ordinal, and minimality gives $\sup B=\nu$.
Since $B$ is a countable set of cardinals,
$\nu$ is a singular cardinal of countable cofinality, and
$\mu<\nu<\delta$.
The restriction of $U$ to $B$ contains every cobounded subset
of $B$. Its dual ideal $J$ therefore extends the bounded ideal,
and
\[
               \tcf(\prod B/J)=\lambda\leq\pp(\nu).
\]
Thus every term in \eqref{eq:bounded} is at most
$\max\{\delta,H_{\mu,\delta}\}$.

For the reverse inequality, let $\nu\in\mathcal S_{\mu,\delta}$,
and let a regular $\lambda$ be represented at $\nu$ by
$\tcf(\prod B/J)=\lambda$. Since $J$ extends the bounded ideal, we can assume $\mu<\min(B)$, and 
therefore $B$ is a countable subset of $\Reg\cap(\mu,\delta)$ bounded by $\nu<\delta$, and $\lambda\in\pcf(B)$. 
Thus, 
$$\lambda\leq\max\pcf(B)\leq\chi(X)=t(X).$$
Taking the supremum over the representations of each $\nu$, and
then over $\nu\in\mathcal S_{\mu,\delta}$, gives the lower bound
$H_{\mu,\delta}$. Together with the lower bound $\delta$ already
obtained, this gives $\chi(X)=t(X)=\max\{\delta,H_{\mu,\delta}\}$.
If $\mathcal S_{\mu,\delta}$ is cofinal in $\delta$, the inequalities
$\pp(\nu)\geq\nu^+$ give $H_{\mu,\delta}\geq\delta$, so the
term $\delta$ may be omitted.
\end{proof}

\begin{proof}[Proof of Corollary~\ref{cor:topological-intro}]
By Theorem~\ref{thm:normality-main}, $X$ is nonnormal precisely when
some nonempty countable $A\subseteq(\mu,\delta)\cap\Reg$ with
$\sup A<\delta$ satisfies $\max\pcf(A)\geq\delta$. Since
$\max\pcf(A)$ is regular and $\delta$ is singular, this inequality
is strict. It follows from~\eqref{eq:bounded} that
\[
    X\text{ is nonnormal}\quad\Longleftrightarrow\quad\chi(X)>\delta.
\]
Theorem~\ref{thm:global} gives $\psi(X)=\delta$ and
$\chi(X)=t(X)\geq\delta$, so the four conditions are equivalent.
\end{proof}

\begin{remark}
If $X$ is nonnormal, choose $A$ as in the preceding proof. The
construction in the proof of Theorem~\ref{thm:global} gives a point
$x\in X$ with $A_x=A$, so by Lemma~\ref{lem:local},
\[
    \psi(x,X)=\sup A<\delta<\max\pcf(A)=\chi(x,X)=t(x,X).
\]
Conversely, any such point forces $\chi(X)>\delta$, so $X$ is nonnormal.
\end{remark}

\end{document}